\documentclass{amsart}
\usepackage{amscd}
\usepackage{amsmath}
\usepackage{latexsym}
\usepackage{verbatim}
\usepackage[latin1]{inputenc}
\usepackage[dvips]{graphics}
\usepackage{graphicx}
\usepackage{amsfonts}
\usepackage{amssymb}
\newtheorem{Thm}{Theorem}[section]
\newtheorem{Prop}[Thm]{Proposition}
\newtheorem{Lem}[Thm]{Lemma}

\theoremstyle{definition}

\theoremstyle{remark}
\newtheorem{Rem}[Thm]{Remark}

\newtheorem{Quest}[Thm]{Question}

\begin{document}
\title[Homological realization]{Homological realization of prescribed\\abelian groups via $K$-theory}
\author{A.\ J. Berrick and M. Matthey}
\address{Department of Mathematics, National University of Singapore, Singapore 117543,
Republic of Singapore} \email{berrick@math.nus.edu.sg}
\address{University of Lausanne, IGAT (Institute for Geometry, Algebra and Topology), BCH,
EPFL, CH-1015 Lausanne, Switzerland}
\email{michel.matthey@ima.unil.ch}
\thanks{Research partially supported by NUS Research Grant R-146-000-049-112 for both
authors, and partially by Swiss National Science Foundation for
the second author} 
\date{28 March, 2004}
\begin{abstract}Using algebraic and topological $K$-theory together with complex $C^{*}$-algebras,
we prove that every abelian group may be realized as the centre of
a strongly torsion generated group whose integral homology is zero
in dimension one and isomorphic to two arbitrarily prescribed
abelian groups in dimensions two and three.
\end{abstract}
\maketitle

%
%
%

%
%
%

%
%
%

%
%
%

%
%
%

%
%
%

\section{Introduction and statement of the main results}

\label{s1}
%
%
%

%
%
%

The main theorem of this paper combines two genres that we now describe.

\smallskip

(1) The first is the ``inverse realization problem'' for functors
taking group-theoretic values. The oldest example, still open,
asks which finite groups can occur as Galois groups of rational
polynomials~\cite{MM}. Another example is the theorem that every
abelian group can be the ideal class group of a Dedekind
domain~\cite{Claborn,Leedham-Green}. Eilenberg and Mac\thinspace
Lane solved the problem for the homotopy group
functors~\cite{EilMacL}. For homology (always integral in this
note), Baumslag, Dyer and Miller showed that every sequence of
abelian groups could be realized as the reduced homology of a
(discrete) group~\cite{BDM}. When one requires the group to be
rich in torsion, the matter becomes more delicate. For finite
groups, for instance, there are well-known constraints due to
Maschke (cohomological version quoted in \cite[p.~227]{HilSta}),
Evens~\cite{Ev}, and Swan~\cite{Swan}. In the same vein, recall
Milgram's counterexample in~\cite{Milgram} to the conjecture
(attributed to Loday) that no nontrivial finite group can have its
first three positive-dimensional homology groups zero,
see~\cite{Giffen}.
 Since any group with a series of finite length whose
factors are either infinite cyclic or locally finite has the
direct sum of all its reduced homology groups either infinite or
zero~\cite{BerrKrop}, it is apparent that one needs to focus on a
more general class of groups with torsion.

\smallskip

Progress has been made with the class of torsion-generated groups, wherein
every element is a product of elements of finite order. In this case, there
is a vestigial version of the previous result\thinspace: if the sequence of
homology groups of a torsion-generated group is finite, then the group
itself is perfect (that is, the first homology vanishes)~\cite{BeMi}.
After a partial result in~\cite{Ber: JAlg}, exploiting the results on
the ideal class group referred to above, the problem was settled
by~\cite{BeMi}, as follows\thinspace:

\smallskip\emph{Let $A_{2},A_{3},\ldots$ be a sequence of abelian groups. Then
there exists a strongly torsion generated group $G$ such that $H_{n}(G)\cong
A_{n}$ for all $n\geq2$.} \smallskip

A \emph{strongly torsion generated} group $G$ is one with the property that,
for each $n\geq2$, there is an element $g_{n}$ of order $n$ that normally
generates $G$, in other words, every element of $G$ is a product of conjugates
of $g_{n}$. The constraint $n\geq2$ in the above statement occurs because such
groups are necessarily perfect \cite[Lem.~7]{BeMi}. Various properties of the
class of strongly torsion generated groups are discussed in~\cite{BeMi}. It
was introduced in~\cite{Ber: JAlg} because its most notable examples arise in
connection with algebraic $K$-theory. They include the infinite alternating
group $A_{\infty}$ and the infinite special linear groups $\operatorname{SL}%
(\mathbb{Z})$ and $\operatorname{SL}(K)$ for any field $K$.
The proof of the above theorem combines techniques of combinatorial group
theory with Miller's affirmation, in~\cite{Miller}, of the Sullivan Conjecture
in homotopy theory.

\smallskip

(2) The second class of results that forms background to the
present work consists of the embedding theorems in combinatorial
group theory. For half a century it has been known that every
group embeds in an algebraically closed group; these groups are
strongly torsion generated~\cite{Neumann,Scott}. When the embedded
group is abelian, it is natural to attempt to embed it as the
centre of the larger group. An embedding as the centre of a
strongly torsion generated group was achieved in~\cite{Ber: JAlg},
again by means of an algebraic $K$-theory use of the result
mentioned earlier on the ideal class group, and in~\cite{BeMi} by
means of combinatorial group theory and homotopy theory. (See
also~\cite{Hickin,Phillips} for constructions when the abelian
group is locally finite.)

\smallskip

Whereas past displays of abelian groups as homology groups have
formed separate results from realizations as centres, here we are
able to combine these two strands in a single realization theorem,
as follows.

\begin{Thm}
\label{principal}
%
%
%
%
%
Let $A$, $B$ and $C$ be three abelian groups. Then, there exists a group $S$
with the following properties\thinspace:

\begin{itemize}
\item [(i)]$S$ is strongly torsion generated;

\item[(ii)] the centre of $S$ is isomorphic to $A$, that is, $\mathcal{Z}%
(S)\cong A$;

\item[(iii)] $S$ is perfect, that is, $H_{1}(S)=0$;

\item[(iv)] the second homology of $S$ is isomorphic to $B$, that is,
$H_{2}(S)\cong B$;

\item[(v)] the third homology of $S$ is isomorphic to $C$, that is,
$H_{3}(S)\cong C$.
\end{itemize}
\end{Thm}

\smallskip

The construction of $S$ is presented in Section~\ref{s4}, followed in
Section~\ref{s5} by the proof of the theorem. In Section~\ref{s6}, we collect
further information on $S$ as a second theorem. The approach is based on a specific
idea of~\cite{Ber: JAlg} -- related to the functors $K_{1}^{\mathrm{a}%
\hspace*{-0.05em}\mathrm{l}\hspace*{-0.05em}\mathrm{g}}$ and $K_{2}%
^{\mathrm{a}\hspace*{-0.05em}\mathrm{l}\hspace*{-0.05em}\mathrm{g}}$ for rings
-- and on some known results on both the topological and algebraic $K$-theory
of complex $C^{\ast}$-algebras; this is all recalled in Sections \ref{s2} and
\ref{s3}. Section~\ref{s7} poses some open questions on the subject, including
two that are planned as the basis for further studies.

\medskip

\noindent\textbf{Acknowledgements\thinspace:} The second named author wishes
to express his deep gratitude to Paul Balmer for fruitful discussions.\smallskip

%
%
%

\section{Recollection on topological $K$-theory of $C^{\ast}$-algebras}

\label{s2}
%
%
%

%
%
%

The present section is devoted to some preparatory material on topological
$K$-theory, needed for the proofs of our main results.

\medskip

In the present paper, by a $C^{\ast}$-algebra, we always mean a \emph{complex}
$C^{\ast}$-algebra. Recall that the topological $K$-theory of $C^{\ast}%
$-algebras has the following properties\thinspace: it is \emph{additive}, that
is, $K_{\ast}^{\mathrm{t}\hspace*{-0.05em}\mathrm{o}\hspace*{-0.05em}%
\mathrm{p}}(\mathcal{A}_{1}\times\mathcal{A}_{2})\cong K_{\ast}^{\mathrm{t}%
\hspace*{-.05em}\mathrm{o}\hspace*{-0.05em}\mathrm{p}}(\mathcal{A}%
_{1})\oplus K_{\ast}^{\mathrm{t}\hspace*{-0.05em}\mathrm{o}\hspace
*{-0.05em}\mathrm{p}}(\mathcal{A}_{2})$; it satisfies \emph{Bott periodicity},
\textsl{i.e.\ }$K_{\ast}^{\mathrm{t}\hspace*{-0.05em}\mathrm{o}\hspace
*{-0.05em}\mathrm{p}}(\mathcal{A})\cong K_{\ast+2}^{\mathrm{t}\hspace
*{-0.05em}\mathrm{o}\hspace*{-0.05em}\mathrm{p}}(\mathcal{A})$; it is
\emph{continuous}, namely, it commutes with filtered colimits (direct limits
of $C^{\ast}$-algebras); and it is \emph{Morita invariant} in the sense that there
is an isomorphism $K_{\ast}^{\mathrm{t}\hspace*{-0.05em}\mathrm{o}%
\hspace*{-0.05em}\mathrm{p}}\big(M_{n}(\mathcal{A})\big)\cong
K_{\ast}^{\mathrm{t}\hspace*{-0.05em}\mathrm{o}\hspace*{-0.05em}\mathrm{p}}(\mathcal{A})$.
Let $\mathcal{K}\cong\mathop{\rm colim}_{n}M_{n}(\mathbb{C})$ be
the $C^{\ast}$-algebra of compact operators on a separable complex
Hilbert space, and $\widehat{\otimes}$ the minimal
(i.e.\ spatial) tensor product of
$C^{\ast}$-algebras. Since $M_{n}(\mathcal{A})\cong\mathcal{A}%
\widehat{\otimes}M_{n}(\mathbb{C})$, combining Morita invariance
and continuity, we deduce that topological $K$-theory is
\emph{stable}, in the
sense that there is an isomorphism $K_{\ast}^{\mathrm{t}\hspace*{-0.05em}%
\mathrm{o}\hspace*{-0.05em}\mathrm{p}}(\mathcal{A}\widehat{\otimes}%
\mathcal{K})\cong K_{\ast}^{\mathrm{t}\hspace*{-0.05em}\mathrm{o}\hspace*
{-.05em}\mathrm{p}}(\mathcal{A})$. The additivity, Bott and stability
isomorphisms are canonical and natural. For $n\in\mathbb{Z}$, we also recall that
\[
K_{n}^{\mathrm{t}\hspace*{-0.05em}\mathrm{o}\hspace*{-0.05em}\mathrm{p}}%
(\mathbb{C})\cong\renewcommand{\arraystretch}{1.2}\arraycolsep1pt\left\{
\begin{array}
[c]{ll}%
\mathbb{Z}\,,\; & \mbox{if $n$ is even}\\
0\,,\; & \mbox{if $n$ is odd}%
\end{array}
\right.  \quad\mbox{so that}\quad K_{n}^{\mathrm{t}\hspace*{-0.05em}%
\mathrm{o}\hspace*{-0.05em}\mathrm{p}}(\mathcal{K})\cong\renewcommand
{\arraystretch}{1.2}\arraycolsep1pt\left\{
\begin{array}
[c]{ll}%
\mathbb{Z}\,,\; & \mbox{if $n$ is even}\\
0\,,\; & \mbox{if $n$ is odd}.
\end{array}
\right.
\]
For more details on $C^{\ast}$-algebras and their topological $K$-theory (in
particular for the properties we have recalled), we refer, for instance, to
the books~\cite{RoLaLa} and~\cite{WeOl}.

\smallskip

We need the following standard result from the theory of $C^{\ast}$-algebras
and their topological $K$-theory. For a proof, we refer to \cite[Ex.\ 9.H,
pp.\ 173--174]{WeOl}, where it is attributed to Higson and Brown; see
also~\cite{RoLaLa} for $M$ countable.

\begin{Prop}
\label{KC}
%
%
%
For any abelian group $M$, there exists a $C^{\ast}$-algebra $\mathcal{E}_{M}%
$, whose topological $K$-theory is given by
\[
K_{2n}^{\mathrm{t}\hspace*{-0.05em}\mathrm{o}\hspace*{-0.05em}\mathrm{p}%
}(\mathcal{E}_{M})\cong M\qquad\mbox{and}\qquad K_{2n+1}^{\mathrm{t}%
\hspace*{-0.05em}\mathrm{o}\hspace*{-0.05em}\mathrm{p}}(\mathcal{E}%
_{M})=0\quad\qquad(n\in\mathbb{Z})\,.
\]
\end{Prop}

\begin{Rem}
\label{product}
%
%
%
By construction, the $C^{\ast}$-algebra $\mathcal{E}_{M}$ of Proposition~\ref{KC} is
\emph{not} unital.
\end{Rem}

\begin{Rem}
\label{functoriality}
%
%
%
The proof we give below of Theorem~\ref{principal} is presented in such a
way that if there is a construction of the $C^{\ast}$-algebra $\mathcal{E}%
_{M}$ in Proposition~\ref{KC} that is functorial in $M$, then the group $S$ in
Theorem~\ref{principal} is also functorial in the abelian groups $A$, $B$ and
$C$ on which it depends, and the homomorphisms occurring in its statement are all
natural. Such a construction of $\mathcal{E}_{M}$ would certainly be of
independent interest.
\end{Rem}

%
%
%

\section{Recollection on algebraic $K$-theory}

\label{s3}
%
%
%

%
%
%

We now proceed with a recollection of standard -- though sometimes highly
nontrivial -- results on algebraic $K$-theory needed in the proof of
Theorem~\ref{principal}, that is presented in Section~\ref{s5} below.
As general references for algebraic $K$-theory, we mention the books
\cite{Ber-Pitman}, \cite{Mag} and \cite{Ros}.

\medskip

Let $\Lambda$ be a \emph{unital} ring. Denote by $\operatorname{GL}%
(\Lambda)=\bigcup_{n\geq1} \operatorname{GL}_{n}(\Lambda)$, $\operatorname{E}%
(\Lambda)=\bigcup_{n\geq1}\operatorname{E}_{n}(\Lambda)$ and
$\operatorname{St}(\Lambda) =\mathop{\rm colim}_{n\geq3}\operatorname{St}%
_{n}(\Lambda)$ the group of infinite invertible matrices, the group of
infinite elementary matrices, and the infinite Steinberg group respectively.
By definition, we have $K_{1}^{\mathrm{a}\hspace*{-.05em}\mathrm{l}%
\hspace*{-.05em}\mathrm{g}}(\Lambda):=\operatorname{GL}(\Lambda
)/\operatorname{E}(\Lambda)$, and, by the Whitehead Lemma, the equalities
$[\operatorname{GL}(\Lambda),\operatorname{GL}(\Lambda)]=[\operatorname{E}%
(\Lambda),\operatorname{E}(\Lambda)]= \operatorname{E}(\Lambda)$ hold (see
Milnor \cite[Lem.\ 3.1]{Mil}); in particular, the group $\operatorname{E}(\Lambda)$ is
perfect and $K_{1}^{\mathrm{a}\hspace*{-.05em}\mathrm{l}\hspace*{-.05em}%
\mathrm{g}}(\Lambda)=H_{1}\big(\operatorname{GL}(\Lambda)\big)$. By definition
of $K_{2}^{\mathrm{a}\hspace*{-.05em}\mathrm{l}\hspace*{-.05em}\mathrm{g}}$
(\textsl{cf.\ }\cite[p.\ 40]{Mil}), we have a functorial exact sequence
\[
0\longrightarrow K_{2}^{\mathrm{a}\hspace*{-.05em}\mathrm{l}\hspace
*{-.05em}\mathrm{g}}(\Lambda)\longrightarrow\operatorname{St}(\Lambda
)\overset{\varphi_{\!_{\Lambda}}}{\longrightarrow}\operatorname{GL}%
(\Lambda)\longrightarrow K_{1}^{\mathrm{a}\hspace*{-.05em}\mathrm{l}%
\hspace*{-.05em}\mathrm{g}}(\Lambda)\longrightarrow0
\]
with $\mathop{\rm Im}\nolimits(\varphi_{\!_{\Lambda}})=\operatorname{E}%
(\Lambda)$. There are isomorphisms $K_{2}^{\mathrm{a}\hspace*{-.05em}%
\mathrm{l}\hspace*{-.05em}\mathrm{g}}(\Lambda)
\cong\mathcal{Z}\big (\operatorname{St}(\Lambda)\big)\cong
H_{2}\big(\operatorname{E}(\Lambda)\big )$, and
$\operatorname{St}(\Lambda)$ has vanishing $H_{1}$ and $H_{2}$,
and is the universal central extension of
$\operatorname{E}(\Lambda)$, see \cite[Thms.\ 5.1 and 5.10]{Mil}.
For a perfect group $P$, one has
$\mathcal{Z}\big(P/\mathcal{Z}(P)\big )=0$ (see for instance
\cite[end of \S\,2.2]{Ber3}). Therefore,
$\operatorname{E}(\Lambda)\cong\operatorname{St}(\Lambda)/\mathcal{Z}\big(
\operatorname{St}(\Lambda)\big)$ is centreless.

\smallskip

It is well-known that $K_{n}^{\mathrm{a}\hspace*{-0.05em}\mathrm{l}%
\hspace*{-0.05em}\mathrm{g}}(\Lambda)$ is isomorphic to $\pi_{n}\big
(B\!\operatorname{GL}(\Lambda)^{+}\big)$ for $n\geq1$ (this is even the
definition for $n\geq3$), to $\pi_{n}\big(B\!\operatorname{E}(\Lambda)^{+}%
\big)$ for $n\geq2$, and to $\pi_{n}\big(B\!\operatorname{St}(\Lambda)^{+}%
\big)$ for $n\geq3$ (see \cite[Cor.\ 5.2.8]{Ros}). Recall from
\cite[Thm.\ 5.2.2]{Ros} that $H_{\ast}(X)\cong H_{\ast}(X^{+})$ holds for
any connected CW-complex $X$, as for example $B\!\operatorname{GL}(\Lambda
)$, $B\!\operatorname{E}(\Lambda)$ and $B\!\operatorname{St}(\Lambda)$. In
particular, knowing that $B\!\operatorname{E}(\Lambda)^{+}$ is $1$-connected
and that $B\!\operatorname{St}(\Lambda)^{+}$ is $2$-connected (\textsl{cf.\ }%
\cite[Thm.\ 5.2.2]{Ros}), by the Hurewicz Theorem \cite[Thm.\ 10.25]{Swi}, the
Hurewicz homomorphism  induces the following epimorphisms and isomorphism\thinspace:
\[
K_{3}^{\mathrm{a}\hspace*{-0.05em}\mathrm{l}\hspace*{-0.05em}\mathrm{g}%
}(\Lambda)\,-\!\!\!\!\!\twoheadrightarrow H_{3}\big(\operatorname{E}%
(\Lambda)\big)\,,\quad\;K_{3}^{\mathrm{a}\hspace*{-0.05em}\mathrm{l}%
\hspace*{-0.05em}\mathrm{g}}(\Lambda)\cong H_{3}\big(\operatorname{St}%
(\Lambda)\big)\quad\;\mbox{and}\quad\;K_{4}^{\mathrm{a}\hspace*{-0.05em}%
\mathrm{l}\hspace*{-0.05em}\mathrm{g}}(\Lambda)\,-\!\!\!\!\!\twoheadrightarrow
H_{4}\big
(\operatorname{St}(\Lambda)\big)%
\]
(the isomorphism is Gersten's Theorem~\cite{Ger}). All indicated
isomorphisms and epimorphisms are canonical and natural. From
\cite[Lem.~1 and proof of Thm.~A]{Ber: JAlg}, we also quote
that\thinspace:

\smallskip\emph{For a unital ring $\Lambda$, the groups $\operatorname{E}(\Lambda)$ and $\operatorname
{St}(\Lambda)$ are strongly torsion generated.}\smallskip

For the definition of negative $K$-theory of a unital ring $\Lambda$,
$K_{-n}^{\mathrm{a}\hspace*{-0.05em}\mathrm{l}\hspace*{-0.05em}\mathrm{g}%
}(\Lambda)$ with $n>0$, we refer to \cite[Def.\ 3.3.1]{Ros}. If $I$ is a
nonunital ring, following \cite[Def.\ 1.5.6]{Ros}, we define the \emph{minimal
unitalization} $\widetilde{I}$ of $I$ as the unital ring given, as a
$\mathbb{Z}$-module, by the direct sum $\widetilde{I}:=I\oplus\mathbb{Z}$, and
equipped with the multiplication given by
\[
(x,\lambda)\cdot(x^{\prime},\lambda^{\prime}):=\big(xx^{\prime}+\lambda
x^{\prime}+\lambda^{\prime}x,\lambda\lambda^{\prime}\big)\,,\quad\mbox
{for}\;\;x,x^{\prime}\in I\;\,\mbox{and}\;\,\lambda,\lambda^{\prime}%
\in\mathbb{Z}\,.
\]
As in \cite[Def.\ 1.5.7]{Ros}, there is a split short exact sequence of
nonunital rings
\[
0\longrightarrow I\longrightarrow\widetilde{I}\overset{\curvearrowleft
}{\longrightarrow}\mathbb{Z}\longrightarrow0\,,
\]
and one defines $K_{\ast}^{\mathrm{a}\hspace*{-0.05em}\mathrm{l}%
\hspace*{-0.05em}\mathrm{g}}(I)$ as to be the kernel of the map $K_{\ast
}^{\mathrm{a}\hspace*{-0.05em}\mathrm{l}\hspace*{-0.05em}\mathrm{g}%
}(\widetilde{I})\longrightarrow K_{\ast}^{\mathrm{a}\hspace*{-0.05em}%
\mathrm{l}\hspace*{-0.05em}\mathrm{g}}(\mathbb{Z})$ induced by the unital ring
homomorphism $\widetilde{I}\longrightarrow\mathbb{Z}$. This construction is
functorial for nonunital ring homomorphisms. For the definition of the
relative $K$-groups $K_{n}^{\mathrm{a}\hspace*{-0.05em}\mathrm{l}%
\hspace*{-0.05em}\mathrm{g}}(\Lambda,J)$, where $J$ is a two-sided ideal in
the unital ring $\Lambda$, we refer, for $n\geq0$, to \cite[Defs.\ 1.5.3 and
5.2.14]{Ros}; for $n>0$, one sets $K_{-n}^{\mathrm{a}\hspace*{-0.05em}%
\mathrm{l}\hspace*{-0.05em}\mathrm{g}}(\Lambda,J):=K_{-n}^{\mathrm{a}%
\hspace*{-0.05em}\mathrm{l}\hspace*{-0.05em}\mathrm{g}}(J)$ (hiding the fact
that $K_{-n}^{\mathrm{a}\hspace*{-0.05em}\mathrm{l}\hspace*{-0.05em}%
\mathrm{g}}$ satisfies excision), see \cite[Def.\ 3.3.1]{Ros}. The above
split exact sequence induces a canonical isomorphism
\[
K_{\ast}^{\mathrm{a}\hspace*{-0.05em}\mathrm{l}\hspace*{-0.05em}\mathrm{g}%
}(\widetilde{I})\cong K_{\ast}^{\mathrm{a}\hspace*{-0.05em}\mathrm{l}%
\hspace*{-0.05em}\mathrm{g}}(\widetilde{I},I)\oplus K_{\ast}^{\mathrm{a}%
\hspace*{-0.05em}\mathrm{l}\hspace*{-0.05em}\mathrm{g}}(\mathbb{Z})\,,
\]
as follows from the long exact sequence in algebraic $K$-theory, see
\cite[Thm.\ 3.3.4]{Ros}. Since $K_{0}^{\mathrm{a}\hspace*{-0.05em}%
\mathrm{l}\hspace*{-0.05em}\mathrm{g}}$ satisfies excision too (see
\cite[Thm.\ 1.5.9]{Ros}), and since \emph{the ring $\mathbb{Z}$ is regular},
so that its negative algebraic $K$-groups all vanish (see
\cite[Ex.\ 3.1.2\,(4) and Def.\ 3.3.1]{Ros}), we get
\[
K_{-n}^{\mathrm{a}\hspace*{-0.05em}\mathrm{l}\hspace*{-0.05em}\mathrm{g}%
}(\widetilde{I})\cong\renewcommand{\arraystretch}{1.2}\arraycolsep1pt\left\{
\begin{array}
[c]{ll}%
K_{-n}^{\mathrm{a}\hspace*{-0.05em}\mathrm{l}\hspace*{-0.05em}\mathrm{g}%
}(I)\,,\; & \mbox{if $n>0$}\\
K_{0}^{\mathrm{a}\hspace*{-0.05em}\mathrm{l}\hspace*{-0.05em}\mathrm{g}%
}(I)\oplus\mathbb{Z}\,,\; & \mbox{if $n=0$}\,.
\end{array}
\right.
\]

\smallskip

Recall that the \emph{cone} of $\mathbb{Z}$ is the unital ring $C(\mathbb{Z})$
consisting of the infinite matrices $(a_{ij})_{i,j\in\mathbb{N}}$ with only
finitely many non-zero (integer-valued) entries in each row and in each
column. The \emph{suspension} of $\mathbb{Z}$ is the quotient $S(\mathbb{Z}%
):=C(\mathbb{Z})/\mathbb{M}(\mathbb{Z})$, where $\mathbb{M}(\mathbb{Z})$ is
the two-sided ideal of finite matrices, \textsl{i.e.\ }the union
$\bigcup_{n\geq1}\mathbb{M}_{n}(\mathbb{Z})$ in $C(\mathbb{Z})$. The main
feature of $C(\mathbb{Z})$ is that it is unital with vanishing algebraic
$K$-theory (including in negative degree). For $k\geq1$, the \emph{$k$-fold
suspension} of a unital ring $\Lambda$ is the unital ring $S^{k}%
(\Lambda):=\Lambda\otimes_{\mathbb{Z}}S(\mathbb{Z})^{\otimes_{\mathbb{Z}}\,k}$
(with the obvious ring structure). The ring $S^{k}(\Lambda)$ satisfies the
following property\thinspace:
\[
K_{n}^{\mathrm{a}\hspace*{-0.05em}\mathrm{l}\hspace*{-0.05em}\mathrm{g}}%
\big(S^{k}(\Lambda)\big)\cong K_{n-k}^{\mathrm{a}\hspace*{-0.05em}\mathrm{l}%
\hspace*{-0.05em}\mathrm{g}}(\Lambda)\qquad(n\in\mathbb{Z})\,.
\]
We also need the fact that algebraic $K$-theory is \emph{additive} in the
sense that there is a natural isomorphism $K_{\ast}^{\mathrm{a}\hspace
*{-0.05em}\mathrm{l}\hspace*{-0.05em}\mathrm{g}}(\Lambda_{1}\times\Lambda
_{2})\cong K_{\ast}^{\mathrm{a}\hspace*{-0.05em}\mathrm{l}\hspace
*{-0.05em}\mathrm{g}}(\Lambda_{1})\oplus K_{\ast}^{\mathrm{a}\hspace
*{-0.05em}\mathrm{l}\hspace*{-0.05em}\mathrm{g}}(\Lambda_{2})$, for any two
unital rings $\Lambda_{1}$ and $\Lambda_{2}$. This property is clear in degree
zero and then follows from the definition in negative degrees; for positive
degrees, see \cite[Prop.\ 1.2.3]{Lod}. One further has a canonical
decomposition $\operatorname{GL}(\Lambda_{1}\times\Lambda_{2})\cong
\operatorname{GL}(\Lambda_{1})\times\operatorname{GL}(\Lambda_{2})$,
and similarly for $\operatorname{E}(-)$ and $\operatorname{St}(-)$, see
\cite[p.\ 326 and Prop.\ 12.8]{Mag}.

\smallskip

We have to discuss $C^{\ast}$-algebras in connection with algebraic
$K$-theory. A $C^{\ast}$-algebra $\mathcal{A}$ is called \emph{stable} if it
is $\ast$-isomorphic to $\mathcal{A}\widehat{\otimes}\mathcal{K}$. Since
$\mathcal{K}\cong\mathcal{K}\widehat{\otimes}\mathcal{K}$, for any
$C^{\ast}$-algebra $\mathcal{A}$, the $C^{\ast}$-algebra $\mathcal{A}%
\widehat{\otimes}\mathcal{K}$ is stable. We now recall a deep result,
namely the \emph{Karoubi Conjecture} (proved in
Suslin-Wodzicki~\cite{SuWo1,SuWo2} -- see Remark~\ref{Rem1} below)\thinspace:

\smallskip\emph{The canonical ``change-of-$K$-theory map'' $K_{\ast}^{\mathrm{a}
\hspace*{-0.05em}\mathrm{l}\hspace*{-0.05em}\mathrm{g}}(\mathcal{A})\longrightarrow
K_{\ast}^{\mathrm{t}\hspace*{-0.05em}\mathrm{o}\hspace*{-0.05em}\mathrm{p}}(\mathcal{A})$
is an isomorphism, for any stable $C^{\ast}$-algebra
$\mathcal{A}$.}\smallskip

(Note that this includes the negative $K$-groups $K_{-n}^{\mathrm{a}\hspace
*{-0.05em}\mathrm{l}\hspace*{-0.05em}\mathrm{g}}$ and $K_{-n}^{\mathrm{t}%
\hspace*{-0.05em}\mathrm{o}\hspace*{-0.05em}\mathrm{p}}$ with $n>0$.)

\begin{Rem}
\label{Rem1}
%
%
%
For the proof of our main results, we will not need the full power of the
Karoubi Conjecture that $K_{n}^{\mathrm{a}\hspace*{-0.05em}\mathrm{l}%
\hspace*{-0.05em}\mathrm{g}}(\mathcal{A}\widehat{\otimes}\mathcal{K})\cong
K_{n}^{\mathrm{t}\hspace*{-0.05em}\mathrm{o}\hspace*{-0.05em}\mathrm{p}%
}(\mathcal{A}\widehat{\otimes}\mathcal{K})$ for any $C^{\ast}$-algebra
$\mathcal{A}$ and any $n\in\mathbb{Z}$ -- which has been proved, as we have
just mentioned. Indeed, in the proofs, we will consider a certain stable
$C^{\ast}$-algebra $\mathcal{F}$ and will only need the values of its
algebraic $K$-theory for $n\leq0$, because of the occurrence of an iterated
suspension (see below). In 1979, Karoubi himself proved in~\cite{Kar1} that
his conjecture is true for $n\leq0$ -- in fact, this motivated the conjecture.
Later, this was shown for $n=1$ in de la Harpe-Skandalis~\cite{dlHS}, and for
$n=2$ in Karoubi~\cite{Kar2} and also in Higson~\cite{Hig}. In the latter it
is also proved that the Karoubi Conjecture holds for Karoubi-Villamajor's
algebraic $K$-theory. Finally, it was Suslin and Wodzicki who established the
conjecture for Quillen's algebraic $K$-theory, in~\cite{SuWo1,SuWo2}. For
related results, the reader may consult \cite{Ina,InaKan,Tap}.
\end{Rem}

%
%
%

\section{Construction of the group $S$ of Theorem \ref{principal}}

\label{s4}
%
%
%

%
%
%

Here we provide the construction of the group $S$ occurring in
Theorem~\ref{principal} and in its complement, namely Theorem~\ref{complement}
below. This will justify the long recollections of Sections \ref{s2} and
\ref{s3}, since then, almost all properties arising in those statements will
become ``automatic'', by the very construction.

\medskip

To begin with, for an abelian group $M$, consider the nonunital $C^{\ast}%
$-algebra
\[
\mathcal{F}_{\!M}:=\mathcal{E}_{M}\widehat{\otimes}\mathcal{K}\,,
\]
where $\mathcal{E}_{M}\ $is as in Proposition \ref{KC}. By the results quoted
in Section~\ref{s2}, we have, for every $n\in\mathbb{Z}$,
\[
K_{2n}^{\mathrm{a}\hspace*{-0.05em}\mathrm{l}\hspace*{-0.05em}\mathrm{g}%
}(\mathcal{F}_{\!M})\cong K_{2n}^{\mathrm{t}\hspace*{-0.05em}\mathrm{o}%
\hspace*{-0.05em}\mathrm{p}}(\mathcal{E}_{M})\cong M\quad\mbox{and}\quad
K_{2n+1}^{\mathrm{a}\hspace*{-0.05em}\mathrm{l}\hspace*{-0.05em}\mathrm{g}%
}(\mathcal{F}_{\!M})\cong K_{2n+1}^{\mathrm{t}\hspace*{-0.05em}\mathrm{o}%
\hspace*{-0.05em}\mathrm{p}}(\mathcal{E}_{M})=0\,\text{.}%
\]

Now let $A$, $B$ and $C$ be three abelian groups, prescribed as in Theorem
\ref{principal}. We require \emph{unital} rings having the appropriate
algebraic $K$-theory in low dimensions. For this purpose, we let
\[
R_{A}:=S^{4}\big(\widetilde{\mathcal{F}_{\!A}}\big)\,,\quad R_{B}:=S^{4}\big(
\widetilde{\mathcal{F}_{\!B}}\big)\quad\mbox{and}\quad R^{C}:=S^{5}\big
(\widetilde{\mathcal{F}_{\!C}}\big)%
\]
be the $4$-fold (resp.\ $5$-fold) algebraic suspension of the minimal
unitalization of the nonunital rings $\mathcal{F}_{\!A}$ and $\mathcal{F}%
_{\!B}$ (resp.\ $\mathcal{F}^{C}$), see Section~\ref{s3}. Assembling
most of the results recalled in Sections \ref{s2} and \ref{s3}, we obtain,
for $n\leq 3$,
\[
K_{n}^{\mathrm{a}\hspace*{-0.05em}\mathrm{l}\hspace*{-0.05em}\mathrm{g}}%
(R_{A})\cong K_{n}^{\mathrm{t}\hspace*{-0.05em}\mathrm{o}\hspace*{-0.05em}
\mathrm{p}}(\mathcal{E}_{A})\cong\renewcommand{\arraystretch}{1.2}
\arraycolsep1pt\left\{
\begin{array}
[c]{ll}%
A\,,\; & \mbox{if $n\leq2$ is even}\\
0\,,\; & \mbox{if $n\leq3$ is odd}%
\end{array}
\right.
\]
and similarly for $B$, while, for $n\leq 4$,
\[
K_{n}^{\mathrm{a}\hspace*{-0.05em}\mathrm{l}\hspace*{-0.05em}\mathrm{g}}%
(R^{C})\cong K_{n-1}^{\mathrm{t}\hspace*{-0.05em}\mathrm{o}\hspace*{-0.05em}
\mathrm{p}}(\mathcal{E}_{C})\cong\renewcommand{\arraystretch}{1.2}
\arraycolsep1pt\left\{
\begin{array}
[c]{ll}%
0\,,\; & \mbox{if $n\leq4$ is even}\\
C\,,\; & \mbox{if $n\leq3$ is odd}.
\end{array}
\right.
\]
We also note that for each $n\in\mathbb{Z}$, the abelian groups $K_{n}%
^{\mathrm{a}\hspace*{-0.05em}\mathrm{l}\hspace*{-0.05em}\mathrm{g}}(R_{A})$
and $K_{n}^{\mathrm{a}\hspace*{-0.05em}\mathrm{l}\hspace*{-0.05em}\mathrm{g}%
}(R_{B})$ (resp.\ $K_{n}^{\mathrm{a}\hspace*{-0.05em}\mathrm{l}%
\hspace*{-0.05em}\mathrm{g}}(R^{C})$) contain a direct summand isomorphic
to $K_{n-4}^{\mathrm{a}\hspace*{-0.05em}\mathrm{l}\hspace*{-0.05em}\mathrm{g}%
}(\mathbb{Z})$ (resp.\ $K_{n-5}^{\mathrm{a}\hspace*{-0.05em}%
\mathrm{l}\hspace*{-0.05em}\mathrm{g}}(\mathbb{Z})$). Finally, for the group
$S$ that we have to construct, we take
\[
S:=\operatorname{St}(R_{A})\times\operatorname{E}(R_{B})\times\operatorname{St}(R^{C})\,.
\]
Note that $K_{n}^{\mathrm{a}\hspace*{-0.05em}\mathrm{l}\hspace*{-0.05em}\mathrm{g}}%
(R^{C})$ being zero, we have $S\cong\operatorname{St}(R_{A})\times\operatorname{E}(R_{B})
\times\operatorname{E}(R^{C})$.

%
%
%

\section{Proof of Theorem \ref{principal}}

\label{s5}
%
%
%

%
%
%

Before the proof, we introduce the following convenient terminology. We call a group $G$ \emph{$n$-perfect} for
some $n\geq1$ if its reduced integral homology vanishes in dimension $\leq n$,
\textsl{i.e.\ }$\widetilde{H}_{j}(G)=0$ for all $j\leq n$; of course,
\emph{$1$-perfect} is the same as \emph{perfect} in the usual sense, while
\emph{$2$-perfect} is often called \emph{superperfect}. In the literature, the terms
\emph{$n$-connected} and \emph{$n$-acyclic} are also to be found.

\medskip

We may now prove, in turn, the statements (i)--(v) of Theorem~\ref{principal}.

\medskip

(i). The proof of Lemma~11 presented in~\cite{BeMi} actually
establishes the following slightly stronger statement than that
asserted there.

\begin{Lem}
Let $H$ be a simple group that, for each $n\geq2$, has a $2$-perfect subgroup
$L_{n}$ possessing an element of order $n$. Suppose that $G$ is a group
containing $H$ in such a way that the normal closure of $H$ in $G$ is $G$
itself. Then every perfect central extension of $G$ is strongly torsion generated.
\end{Lem}

As in~\cite{BeMi}, this result may be applied to the case where
$G=\operatorname{E}(\Lambda)$ for any unital ring $\Lambda$ (with $H\cong
A_{\infty}$, see {\sl loc.\ cit.}), to yield that every perfect
central extension of $\operatorname{E}(\Lambda)$ is strongly torsion generated.
For the present circumstance, we take $\Lambda=R_{A}\times
R_{B}\times R^{C}$, and deduce that the perfect central extension
$S$ is strongly torsion generated.

\smallskip

(ii). Since for any unital ring $\Lambda$ the group $\operatorname{E}(\Lambda)$
has trivial
centre, while the centre of $\operatorname{St}(\Lambda)$ is precisely $K_{2}%
^{\mathrm{a}\hspace*{-0.05em}
\mathrm{l}\hspace*{-0.05em}\mathrm{g}}(\Lambda )$, we immediately
have $A$ as the centre of $S$.

\smallskip

(iii). As quoted in Section~\ref{s2}, the groups $\operatorname{E}(\Lambda)$ and
$\operatorname{St}(\Lambda)$ are perfect, for any unital ring
$\Lambda$; and a finite product of perfect groups is perfect.
(Here, one can also recall from \cite[Lem.~7]{BeMi} that
\emph{every strongly torsion generated group is perfect}.)

\smallskip

(iv)\,\&\,(v). Lastly, for the claims about the homology groups of
$S$, we observe from Hurewicz isomorphisms that the first nonzero
reduced homology groups of $\operatorname{St}(R_{A})$, of $\operatorname{E}(R_{B})$
and of $\operatorname{St}(R^{C})$ occur in dimensions $4$, $2$ and $3$
respectively. Moreover, combining with the epimorphism in the next
dimension, we have

\begin{itemize}
\item[\raisebox{.2em}{${\scriptscriptstyle\bullet}\!$}] $H_{2}\big(%
\operatorname{E}(R_{B})\big)\cong
K_{2}^{\mathrm{a}\hspace*{-0.05em}\mathrm{l}\hspace
*{-0.05em}\mathrm{g}}(R_{B})\cong B$
\quad and\quad
$H_{3}\big(\operatorname{E}(R_{B})\big)=0$

\item[\raisebox{.2em}{${\scriptscriptstyle\bullet}\!$}] $H_{3}%
\big(\operatorname{St}(R^{C})\big)\cong K_{3}^{\mathrm{a}\hspace*{-0.05em}%
\mathrm{l}\hspace*{-0.05em}\mathrm{g}}(R^{C})\cong C$
\quad and\quad
$H_{4}\big(\operatorname{St}(R^{C})\big)=0$\,.
\end{itemize}

Hence, the desired results are immediate from the K\"{u}nneth
Theorem. This completes the proof. \hfill$\Box\smallskip$

\begin{Rem}
\label{rem}
%
%
%

\begin{itemize}
\item [(i)]Observe that for $n_{\!A}\geq2$ as large as we like, we can replace
$S^{4}(\widetilde{\mathcal{F}_{\!A}})$ by $S^{2n_{\!A}}(\widetilde
{\mathcal{F}_{\!A}})$. Note also that for
$n_{\!A}\geq 3$,
\[
\qquad\quad\;\; H_{4}\big(\operatorname{St}(R_{A})\big)\cong K_{4}^{\mathrm{a}
\hspace*{-0.05em}\mathrm{l}\hspace*{-0.05em}\mathrm{g}}(R_{A})\cong A
\]
(however, for $n_{\!A}=2$ we get $A\oplus\mathbb{Z}$). We can play
the same game with $B$; while for $C$ we can take
$S^{2n_{C}+1}(\widetilde{\mathcal{F}_{\!C}})$ with $n_{C}\geq2$.
We need $n_{\!A},n_{\!B}\geq3$ in
Theorem~\ref{complement}\thinspace(i) below.

\item[(ii)] For $n\geq2n_{\!A}$, $K_{n}^{\mathrm{a}\hspace*{-0.05em}%
\mathrm{l}\hspace*{-0.05em}\mathrm{g}}(R_{A})$ contains a direct summand
$K_{n-2n_{\!A}}^{\mathrm{a}\hspace*{-0.05em}\mathrm{l}\hspace*{-0.05em}%
\mathrm{g}}(\mathbb{Z})$, and so, when $n-2n_{\!A}\equiv 1\
(\mathrm{\operatorname{mod}\;}4)$ and $n-2n_{\!A}\geq5$, a direct
summand isomorphic to $\mathbb{Z}$, see \cite[Thms 5.3.12 and
5.3.13]{Ros}; similarly with $B$ and $C$.

\item[(iii)] When $A=0$, the unital ring $R_{A}$ is isomorphic to the
$2n_{\!A}$-fold suspension of the ring of integers, $S^{2n_{\!A}}(\mathbb{Z}%
)$; and similarly for $R_{B}$ (resp.\ $R^{C}$) when $B=0$ (resp.\ $C=0$).

\item[(iv)] There are two drawbacks to our construction of the group $S$
above\thinspace: first, our construction of $S$ is not functorial in $A$, $B$ and $C$
(see however Remark~\ref{functoriality}); secondly, for $A$, $B$ and $C$
countable (and even finite), the group $S$ is not countable.
\end{itemize}
\end{Rem}

%
%
%

\section{Further consequences of the construction}

\label{s6}
%
%
%

%
%
%

We now show that the group $S$ of Theorem~\ref{principal} can be constructed
in such a way that further properties hold, that are stated as Theorem~\ref{complement}
below.

\medskip

In this section, assuming that we have taken
$n_{\!A},n_{\!B}\geq3$ in the notation of
Remark~\ref{rem}\thinspace(i), we prove the next result, which
complements Theorem~\ref{principal}.

\begin{Thm}
\label{complement}
%
%
%
The group $S$ of Theorem~\ref{principal} has the following further
properties\thinspace:

\begin{itemize}
\item [(i)]when $B=0$, one has $H_{4}(S)\cong A$;

\item[(ii)] for infinitely many dimensions $n$, the homology group $H_{n}(S)$
contains an infinite cyclic direct summand;

\item[(iii)] every finite-dimensional complex representation of $S$ is trivial;

\item[(iv)] $S/A$ is strongly torsion generated, centreless, and
has $H_{2}(S/A)\cong A\times B$.
\end{itemize}
\end{Thm}

\noindent\text{Proof.\quad}  As in Remark~\ref{rem}\thinspace(i), we choose
integers $n_{\!A},n_{\!B}\geq3$ and $n_{C}\geq2$, and take
\[
R_{A}=S^{2n_{\!A}}(\widetilde{\mathcal{F}_{\!A}})\,,\quad R_{B}=S^{2n_{\!B}%
}(\widetilde{\mathcal{F}_{\!B}})\quad\mbox{and}\quad R^{C}=S^{2n_{C}+1}(\widetilde
{\mathcal{F}_{\!C}})\,.%
\]

\smallskip

(i). Since $B=0$, by the Hurewicz isomorphism, the first nonzero reduced homology group
of $\operatorname{E}(R_{B})$ is in dimension $2n_{\!B}\geq 6$, so that $H_{4}\big(\operatorname{E}(R_{B})\big)=0$.
(Of course, one could also modify $S$ by omitting all usage of $R_{B}$ in this case.) As
a consequence and since $n_{\!A}\geq 3$, the isomorphism $H_{4}\big(\operatorname{St}(R_{A})
\big)\cong A$ of Remark~\ref{rem}\,(i) combines with the K\"{u}nneth Theorem to give
the result.

\smallskip

(ii). By \cite[Thm.~2.1]{Arl1}, for $n\geq 3$, there are homomorphisms
\[
K_{n}^{\mathrm{a}\hspace*{-0.05em}\mathrm{l}\hspace*{-0.05em}\mathrm{g}%
}(\Lambda)\longrightarrow H_{n}\big(\operatorname{E}(\Lambda)\big)
\longrightarrow K_{n}^{\mathrm{a}\hspace*{-0.05em}\mathrm{l}\hspace*{-0.05em}
\mathrm{g}}(\Lambda)
\;\;\mbox{and}\;\;
K_{n}^{\mathrm{a}\hspace*{-0.05em}\mathrm{l}\hspace*{-0.05em}\mathrm{g}%
}(\Lambda)\longrightarrow H_{n}\big(\operatorname{St}(\Lambda)\big)\longrightarrow
K_{n}^{\mathrm{a}\hspace*{-0.05em}%
\mathrm{l}\hspace*{-0.05em}\mathrm{g}}(\Lambda)\,,
\]
such that in each case the composite is multiplication by a positive integer.
This means that elements of infinite order in $K_{n}^{\mathrm{a}%
\hspace*{-0.05em}\mathrm{l}\hspace*{-0.05em}\mathrm{g}}(\Lambda)$ are mapped
to elements of infinite order by the composite, and so must have images of
infinite order in $H_{n}\big(\operatorname{E}(\Lambda)\big)$ and in $H_{n}\big(\operatorname{St}%
(\Lambda)\big)$, respectively. Hence, by
Remark~\ref{rem}\thinspace(ii), no matter which of the possible
values for $n_{\!A}$, $n_{\!B}$ and $n_{C}$ we choose in our
construction, this yields an infinite cyclic direct summand (and
possibly three such summands) in $H_{n}\big(\operatorname{E}(\Lambda)\big)$, in
$H_{n}\big (\operatorname{St}(\Lambda)\big)$, and therefore in $H_{n}(S)$, for
infinitely many values of $n$.

\smallskip

(iii). Because a surjective unital ring homomorphism, such as that from the cone of a
ring to its suspension, induces a surjection of Steinberg groups,
it follows from the construction that the group $S$ is a
homomorphic image of the Steinberg group
\[
\operatorname{St}\Big(C\big(\,S^{2n_{\!A}-1}(\widetilde{\mathcal{F}_{\!A}}%
)\,\times\,S^{2n_{\!B}-1}(\widetilde{\mathcal{F}_{\!B}})\,\times\,S^{2n_{C}%
}(\widetilde{\mathcal{F}_{\!C}})\,\big)\Big)\,.%
\]
This group is acyclic and torsion-generated; therefore, by the main result
of~\cite{Ber: BLMS}, it has no nontrivial finite-dimensional complex
representation (\textsl{cf.\ }\cite[p.~191]{Ber: JAlg}).

\smallskip

(iv). We have
$\operatorname{St}(R_{A})/A\cong\operatorname{E}(R_{A})$ whence it
follows that the quotient $S/A$ is isomorphic to
$\operatorname{E}(R_{A}\times R_{B}\times R^{C})$ and is strongly
torsion generated, centreless, with second homology group
isomorphic to
$K_{2}^{\mathrm{a}\hspace*{-0.05em}\mathrm{l}\hspace*{-0.05em}\mathrm{g}}(R_{A}
\times R_{B}\times R^{C})\cong A\times B$.\hfill$\Box\smallskip$

\smallskip

The interest of item (ii) is heightened by the following
observation. In~\cite{BeMi}, to construct a strongly torsion
generated group $S$ with prescribed centre $A$, one starts with a
centreless strongly generated group $S^{\prime}$ with its reduced
integral homology concentrated in dimension $2$ and isomorphic to
$A$, and then one takes for $S$ the universal central extension of
$S^{\prime}$ (see details in \cite[Proof of Cor.~16]{BeMi}). In
that sense, if $A$ has very few nonvanishing (resp. nontorsion)
integral homology groups (for example, if $A$ is free abelian of
finite rank), then by the Lyndon-Hochschild-Serre spectral
sequence one can expect $S$ to have very few nonvanishing (resp.
nontorsion) integral homology groups as well. The argument
of~\cite{BeMi} obliges one to be in this situation in order to
establish that $S$ is strongly torsion generated. In contrast, for
the present construction of $S$, infinite higher homology groups
are inescapable.

%
%
%

\section{Some open questions}

\label{s7}
%
%
%

%
%
%

The following questions are prompted by Theorem~\ref{principal}. For each of
them, the further requirement of a functorial construction is also of interest.

\begin{Quest}
\label{Q2}
%
%
%
Can one find a \emph{countable} group $S$ as in the statement of
Theorem~\ref{principal} for $A$, $B$ and $C$ finite or countable\thinspace? What
about if we replace, everywhere, the word ``countable'' by ``finitely
generated'', or ``finitely presented''\thinspace?
\end{Quest}

As a matter of comparison, it is known that any group (resp.\ countable group,
finitely generated group, finitely presented group, geometrically finite
group) embeds in an acyclic group (resp.\ countable acyclic group,
seven-generator acyclic group, finitely presented acyclic group, geometrically
finite acyclic group), see~\cite{BDH}. (A group $G$ is called
\emph{geometrically finite} if there exists a model for its classifying space
that is a finite $CW$-complex; in particular, $G$ is then torsion-free.) Since
a strongly torsion generated group, by its very definition, contains ``a lot
of torsion'', Question~\ref{Q2} for ``geometrically finite'' has a negative answer.

\begin{Quest}
\label{Q3}
%
%
%
Given $n\geq3$, is it possible to define an $n$-perfect strongly
torsion generated group $S$ with centre and $H_{n+1}(S)$
prescribed abelian groups\thinspace?
\end{Quest}

The cases $n=1$ and $n=2$ have been achieved above. In the case $n=2$, the extension
$A\rightarrowtail S\,-\!\!\!\!\!\twoheadrightarrow S/A$ is the universal central extension
of the perfect group $S/A$, and so, $H_{2}(S/A)\cong A$.

\smallskip

The case $n=\infty$ (so to speak) of the above question merits special attention.

\begin{Quest}
\label{Q4}
%
%
%
Given an abelian group $A$, is it possible to construct an acyclic
strongly torsion generated group $S$ with centre isomorphic to
$A$\thinspace?
\end{Quest}

Again, we can compare this with the fact, proved in~\cite{Ber2},
that any abelian group is the centre of some acyclic group, in a
functorial and explicit way.

\begin{Quest}
\label{Q5}
%
%
%
Given an abelian group $M$ and $n\geq3$, does there exist a group
$G$ with $G\big/[G,G]\cong M$, and such that $[G,G]$ is
$n$-perfect strongly torsion generated with trivial
centre\thinspace? Can one further require $[G,G]$ to be
acyclic\thinspace?
\end{Quest}

Again, we have already dealt with the case $n=2$, where, as group $G$, we take
the infinite general linear group $\operatorname{GL}\big(S^{3}(\widetilde
{\mathcal{F}_{\!M}})\big)$, which satisfies
\[
G\big/[G,G]=K_{1}^{\mathrm{a}\hspace*{-0.05em}\mathrm{l}\hspace*{-0.05em}%
\mathrm{g}}\big(S^{3}(\widetilde{\mathcal{F}_{\!M}})\big)\cong
M\quad \mbox{and}\quad\lbrack
G,G]=\operatorname{E}\big(S^{3}(\widetilde{\mathcal{F}_{\!M}})\big
)=\operatorname{St}\big(S^{3}(\widetilde{\mathcal{F}_{\!M}})\big)\,.
\]

\smallskip

Finally, here are two questions that we plan to address in subsequent work.

\begin{Quest}
\label{Q: real}
%
%
%
Is it possible similarly to construct a strongly torsion generated
group $S$ with given centre and having more prescribed homology
than in Theorem~\ref{principal}\thinspace?
\end{Quest}

\begin{Quest}
\label{Q: M ring and groups}
%
%
%
What kind of information do the ring $\widetilde{\mathcal{F}_{\!M}}$ and its $K$-theory
convey about the abelian group $M$\thinspace?
\end{Quest}

%
%
%

%
%
%
\end{document}